\documentclass[10pt, a4paper]{article}

\usepackage{amssymb,amsmath,amsthm,amscd,enumerate}

\usepackage[all]{xy}

\usepackage[french]{babel}

\newcommand{\lr}{\longrightarrow}
\newcommand{\T}{{\mathrm Tr}}

\advance\oddsidemargin -1,5cm
\advance\evensidemargin -2cm

\begin{document}
\begin{center}
{\bf On the non additivity of the trace in derived categories}\\

Daniel Ferrand
\end{center}
\vspace{5mm}

In this note  we provide an example of an endomorphism of a short exact sequence of perfect complexes, with the trace of the middle map not equal to the sum of the traces of the two other ones. The point is that the squares involved are commutative only up to homotopy. In view of this example I have found in 1968, Deligne immediately created his "categories spectrales",  and soon afterwards Illusie introduced the "filtered derived categories" where a satisfactory kind of additivity is restored for the trace. 
\vspace{7mm}

{\bf 1. Rappels}\\

R\'eduisons les rappels au minimum. Soit $A$ un anneau commutatif. On ne
consid\`ere que des complexes born\'es de $A$-modules projectifs de type
fini
$$
K = \qquad \cdots \lr \; K^{n-1}\; \stackrel{d}{\lr} \; K^n\;
\stackrel{d}{\lr} \; K^{n+1}\; \lr  \cdots
$$
Dans la suite, ils seront nomm\'es \emph{complexes parfaits} (En toute rigueur,
il faudrait dire : strictement parfaits). Pour un endomorphisme de complexe
parfait
$u : K
\rightarrow K$, on pose
$$
\T(u) \; = \; \sum (-1)^{i}\T(u^{i}).
$$
Consid\'erons trois complexes parfaits, $K, L$ et $M$ et une suite exacte
courte
$0\lr K
\stackrel{j}{\lr }L \stackrel{q}{\lr} M \lr 0$ (ce qui veut dire qu'en chaque degr\'e $n$, la suite $0 \rightarrow K^n \rightarrow L^n \rightarrow M^n \rightarrow 0$ est exacte).
% $j$ et $q$
%induisent en chaque degr\'e $n$ une suite exacte de $A$-modules).
%$$
%0\; \lr \;K^n\; \stackrel{j^n}{\lr} \; L^n\; \stackrel{q^n}{\lr}\; M^n\;
%\lr \; 0.)
%$$
Consid\'erons un endomorphisme de cette suite exacte, c'est-\`a-dire trois
endomorphismes de complexes
$u, v$ et
$w$ tels que le diagramme
$$
\leqno{(1)} \hspace{5,5cm}
\begin{CD}
K @>j>> L @>q>>M\\
@VuVV @VvVV  @VVwV\\
K@>>j> L @>>q> M
\end{CD}
$$
soit \emph{commutatif}. Alors on a la formule d'additivit\'e
$$
\leqno{(2)} \hspace{5,5cm} \T(v) \; =\; \T(u) + \T(w),
$$
puisqu'elle est vraie degr\'e par degr\'e.\\

On v\'erifie sans peine que deux endomorphismes homotopes d'un complexe
parfait ont
la m\^eme trace.  On a pu esp\'erer que la propri\'et\'e d'additivit\'e (2)
serait conserv\'ee lorsque le diagramme (1) ne commute \emph{qu'\`a homotopie
pr\`es}. Il n'en est rien.
\vspace{5mm}

{\bf 2.- L'exemple}\\

On suppose que l'anneau $A$ contient un \'el\'ement $\varepsilon$, non nul
et de
carr\'e nul, et on consid\`ere les complexes (parfaits) suivants :
\begin{itemize}
\item $K$ est $A$ plac\'e en degr\'e $1$ ;
\item $L = (A\stackrel{\varepsilon}{\lr} A)$, en degr\'e $0$  et $1$ ;
\item $M$ est $A$ plac\'e en degr\'e $0$.
\end{itemize}

On a une suite exacte courte de complexes
$$
\begin{CD}
K@>>>L@>>> M\\
&&&&&\\
&& A@=A\\
&&@V{\varepsilon}VV &&\\
A @= A &&
\end{CD}
$$
Finalement, on consid\`ere les trois endomorphismes $u = 0,\; w = 0$, et $v =
{0\choose
\varepsilon} : L\rightarrow L\;$ (c'est-\`a-dire $0$ en degr\'e $0$, et
$\varepsilon$ en degr\'e $1$).  

L'additivit\'e des traces n'est pas respect\'ee puisque $\T(u) = \T(w) =
0$, et
$\T(v) = - \varepsilon$. Or, dans le diagramme
$$
\begin{CD}
K @>>> L @>>>M\\
@V0VV @VvVV  @VV0V\\
K@>>> L @>>> M
\end{CD}
$$
le carr\'e de droite est commutatif, et le carr\'e de gauche est commutatif
\emph{\`a homotopie pr\`es}, comme on le voit en prenant pour homotopie  $h$
l'application identique
$K^1 = A
\lr L^0 = A$, qui est en tirets sur la figure suivante qui r\'esume la
situation :
$$
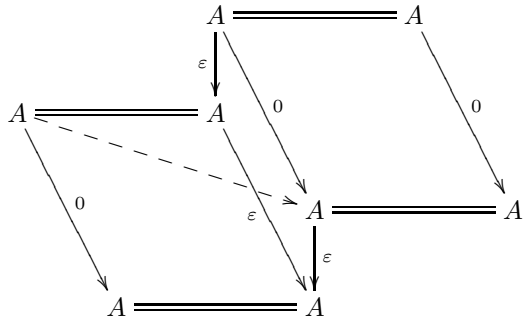
\xymatrix{
& &  A\ar@{=}[rr] \ar[d]_{\varepsilon} \ar[ddr]^0 & & A
\ar[ddr]^0 &\\
A\ar@{=}[rr] \ar@{-->}[drrr] \ar[ddr]^0 & &
A\ar[ddr]_{\varepsilon}  & & &\\
 & & & A \ar[d]^{\varepsilon} \ar@{=}[rr]& & A\\
& A \ar@{=}[rr] & & A & &
}
$$

{\bf 3.- Commentaires}\\

Au milieu des ann\'ees 60 je suivais le S\'eminaire de G\'eom\'etrie
Alg\'ebrique pilot\'e par Grothendieck. Au d\'ebut de 1967 il me
demanda de faire un expos\'e, et de le r\'ediger, sur le d\'eterminant
des complexes parfaits (l'expos\'e XI, manquant dans SGA 6). Ce devait
\^etre, me dit-il, un simple travail de v\'erification \`a partir de ses
notes, lesquelles annon\c{c}aient, entre autre, la multiplicativit\'e du
d\'eterminant pour les morphismes de triangles distingu\'es \emph{ dans
la cat\'egorie d\'eriv\'ee des complexes parfaits} (cette
multiplicativit\'e antra\^ine l'additivit\'e de la trace).

Deux difficult\'es m'emp\^ech\`erent d'avancer.

- Tenir compte des signes - omnipr\'esents -  dont on doit affecter les
isomorphismes li\'es \`a des permutations des facteurs.

- Passer \`a la cat\'egorie d\'eriv\'ee, c'est-\`a-dire essentiellement
envisager des diagrammes qui ne commutent qu'\`a homotopie
pr\`es.\\

Pour surmonter la premi\`ere difficult\'e et rendre compr\'ehensible le
destin de
ces signes, Grothendieck eut, peu apr\`es, l'id\'ee lumineuse de faire du
d\'eterminant, non un simple module inversible, mais un module inversible
\emph{gradu\'e} (concentr\'e en un seul degr\'e si Spec($A$) est connexe) : le
d\'eterminant du complexe parfait
$\cdots
\lr K^n \lr K^{n+1} \lr \cdots $ o\`u $K^n$ est projectif de rang
$r(n)$, devient le couple $(L, r)$, o\`u
$$L = \bigotimes_n
(\wedge^{r(n)}K^n)^{(-1)^n}, \quad \textrm{et}\quad r = \sum
(-1)^nr(n).$$

L'isomorphisme de commutativit\'e $(L, r)\otimes(M,s) \simeq
(M,s)\otimes(L,r)$ doit \^etre affect\'e du signe \og de
Koszul\fg\; $(-1)^{rs}$. Voir [KM] (il s'agit bien du produit $rs$, et
non de la somme, comme une coquille p.20 de cet article pourrait le laisser
croire).\\

Que les questions d'homotopie me g\^enassent \`a ce point montrait,
aux yeux de Grothendieck, mon inaptitude flagrante. Mais, finalement, j'ai
os\'e
penser que c'\'etait peut-\^etre faux, ce qui me conduisit alors tr\`es vite au
contre-exemple trivial signal\'e plus haut. C'\'etait en juillet 1968.
Deligne se
mit au travail, et d\`es septembre, il pr\'esenta ses \og cat\'egories
spectrales\fg \; o\`u les triangles, les \og vrais\fg\,, sont assuj\'ettis
\`a des
conditions fortes  qui rendent une th\'eorie du d\'eterminant possible.

\`A ma connaissance, ce texte n'a pas \'et\'e publi\'e.

Voici un extrait de son introduction : {\it L'id\'ee est la suivante : les
triangles qu'on rencontre \og en pratique\fg\; sont toujours d\'eduits d'un
objet
plus fin : \og un vrai triangle\fg\; et, plus important, les
endomorphismes de triangles qu'on rencontre \og en pratique\fg\;
proviennent toujours d'endomorphismes de vrais triangles. Pour de tels
endomorphismes, la formule [d'additivit\'e des traces] redevient correcte.

Malheureusement, les \og vrais triangles\fg\; de $D(\mathcal{A})$ forment
une cat\'egorie triangul\'ee, de m\^eme que les cat\'egories des \og vrais
triangles\fg\; de \og vrais triangles\fg \ldots Cette machine infernale
explique la complication du formalisme de compatibilit\'e requis pour
d\'efinir les \og cat\'egories spectrales\fg \ldots }\\

Pour mener \`a bien cette r\'edaction il m'aurait fallu alors tout reprendre dans ce nouveau
contexte, et d'abord le comprendre. Une certaine lassitude, et un int\'er\^et pour des
sujets plus directement g\'eom\'etriques m'ont fait abandonner ce travail, et
SGA 6 est paru sans l'expos\'e XI.\\

Peu apr\`es, Illusie introduisit ses \og cat\'egories d\'eriv\'ees
filtr\'ees\fg \;
[Il] ch. V, o\`u {\it \og la construction du c\^one, non fonctorielle, est
remplac\'ee par celle de gradu\'e associ\'e, qui l'est.\fg} Dans ce cadre
l'additivit\'e de la trace est restaur\'ee sous la forme suivante (V 3.7.7,
p.310)
: la trace d'un endomorphisme d'un complexe filtr\'e dont le gradu\'e est
parfait,
est
\'egale \`a la trace de l'endomorphisme induit sur le complexe gradu\'e
associ\'e.\\

En 1975, Knudsen et Mumford publi\`erent un article sur le d\'eterminant des
complexes parfaits [KM]. Ils \'etablissent sa multiplicativit\'e pour les
endomorphismes de suites exactes courtes (p.41), et ils soulignent en note que
les carr\'es en jeu doivent \^etre \og vraiment\fg\; commmutatifs, et non pas
seulement \`a homotopie pr\`es. L'exemple signal\'e plus haut justifie cette mise en garde.
De plus, p.44, ils montrent que si le sch\'ema de
base est \emph{r\'eduit}, il suffit de supposer la commutativit\'e \`a
homotopie
pr\`es. L'emploi d'un nilpotent dans l'exemple est donc in\'evitable.\\

{\bf Bibliographie}\\

[D]\; DELIGNE P., {\it Cat\'egories spectrales}, Manuscrit, (Sept. 1968).

[Il]\; ILLUSIE L.,{\it Complexe cotangent et d\'eformations, I}, LNM 239,
Springer
(1971).

[KM]\; KNUDSEN F. and MUMFORD D., {\it The projectivity of the moduli space of
stable curves.\, I : Preliminaries on "det" and "Div"}, Math. Scand. 39 (1976)
19-55.
\vspace{5mm}

{\sc IRMAR, Universit\' e  de Rennes 1, Campus de Beaulieu, F-35040 Rennes Cedex}

{\it E-mail address : }\texttt {daniel.ferrand}{@}\texttt{univ-rennes1.fr}

\end{document}